\documentclass[a4paper,8pt]{article}
\usepackage{amsmath,amssymb,bbm,theorem}

\newcommand{\assign}{:=}
\newcommand{\noplus}{}
\newcommand{\tmop}[1]{\ensuremath{\operatorname{#1}}}
\newcommand{\tmstrong}[1]{\textbf{#1}}
{\theorembodyfont{\rmfamily}\newtheorem{remark}{Remark}}

\begin{document}

\title{Concavity of Perelman's $\mathcal{W}$-functional over the space of
K\"{a}hler potentials}\author{\\
{\tmstrong{NEFTON PALI}}}\date{}\maketitle

\begin{abstract}
  In this short note we observe that the concavity of Perelman's
  $\mathcal{W}$-functional over a neighborhood of a K\"{a}hler-Ricci soliton inside
  the space of K\"{a}hler potentials is a direct consequence of author's solution
  of the variational stability problem for K\"{a}hler-Ricci solitons.
  Independently, we provide a rather simple proof of this fact based on some
  elementary formulas obtained in our previous work.
\end{abstract}

Let $\left( X, J, \omega \right)$ be a compact K\"{a}hler-Ricci soliton and let
$\omega_t \assign \omega + i \partial \overline{\partial} \varphi_t$ be a family of K\"{a}hler metrics with
$\varphi_0 = 0$. We denote by $f_t$ the unique function such that $i \partial
\overline{\partial} f_t = \omega_t - \tmop{Ric} \left( \omega_t \right)$ and
$\int_X e^{- f_t} \omega^n_t = n!$. We consider the Riemannian metric $g_t
\assign - \omega_t J$ and the positive volume form $\Omega_t \assign e^{- f_t}
\omega^n_t / n!$. We observe that Perelman's $\mathcal{W}$-functional \cite{Per}
satisfies the identity $\mathcal{W} \left( g_t, \Omega_t \right) = 2 \int_X
f_t \Omega_t$. The monotony of this quantity along the K\"{a}hler-Ricci flow was discovered in \cite{Pal1}. The inequality
\begin{equation}
  \label{concav}  \frac{d^2}{d t^2} \vphantom{dt}_{|_{t = 0}} \mathcal{W} (g_t,
  \Omega_t) \leqslant 0,
\end{equation}
follows immediately form the existence of an identification map $\omega_t
\longleftrightarrow J_t\equiv$ complex structure compatible with $\omega$ (see lemma 30 in \cite{Pal2}), from
the diffeomorphism invariance of $\mathcal{W}$ and from author's solution of the
variational stability problem for K\"{a}hler-Ricci solitons \cite{Pal2,Pal3}. These facts imply
also that the equality hold in (\ref{concav}) if and only if $\dot{\varphi}_0
= \psi + \bar{\psi}$ with $\left( \nabla_{g_0} \psi \right)_J^{1, 0}$ 
holomorphic vector field (See proposition 2 and lemma 22 in \cite{Pal2}).

We provide now an independent proof of the inequality (\ref{concav}) and of
the identification of the kernel of the left hand side. The first variation of
the quantity $\mathcal{W}_t \assign \mathcal{W} (g_t, \Omega_t)$ follows form
a computation quite similar to one given in section 4 of \cite{Pal1}. We
include it here for readers convenience. We introduce the $\Omega$-divergence operator acting on
vector fields $\xi$ as
\begin{eqnarray*}
  \tmop{div}^{\Omega} \xi & \assign & \frac{d (\xi \neg \Omega)}{\Omega} .
\end{eqnarray*}
We define the weighted real Laplacian $\Delta^{\Omega}_g u \assign -
\tmop{div}^{\Omega} \nabla_g u$, acting on functions $u$. We can assume without loss of
generality that the potential $\varphi_t$ is normalized in a way that $\int_X
\dot{\varphi}_t \Omega_t \equiv 0$. (Indeed we can replace the potential $\varphi_t$ with
$\tilde{\varphi}_t = \varphi_t - \int^t_0 d s \int_X \dot{\varphi}_s
\Omega_s$). We consider now the
function
\begin{equation}
  \label{volVar}  \dot{\Omega}^{\ast}_t \assign \dot{\Omega}_t / \Omega_t = -
  \frac{1}{2} \Delta_{g_t}  \dot{\varphi}_t - \dot{f}_t .
\end{equation}
Time deriving the normalizing condition $\int_X \Omega_t \equiv 1$ we obtain
\begin{equation}
  \label{volNor} \int_X \dot{\Omega}^{\ast}_t \Omega_t \equiv 0 .
\end{equation}
Time deriving the identity $\omega_t = \tmop{Ric} \left( \Omega_t \right)$ we
obtain $i \partial \overline{\partial}  \dot{\varphi}_t = - i \partial
\overline{\partial} \dot{\Omega}^{\ast}_t$. This combined with the
normalization of $\varphi_t$ and with (\ref{volNor}) implies
\begin{equation}
  \label{expVarVol}  \dot{\varphi}_t = - \dot{\Omega}^{\ast}_t .
\end{equation}
Using (\ref{expVarVol}) and (\ref{volVar}) we obtain
\begin{equation}
  \label{evolf} 2 \dot{f}_t = - \Delta_{g_t}  \dot{\varphi}_t \noplus + 2
  \dot{\varphi}_t .
\end{equation}
We deduce the identities
\begin{eqnarray*}
  \dot{\mathcal{W}}_t & = & 2 \int_X \dot{f}_t \Omega_t + 2 \int_X f_t 
  \dot{\Omega}^{\ast}_t \Omega_t\\
  &  & \\
  & = & - \int_X \Delta_{g_t}  \dot{\varphi}_t \noplus - 2 \int_X f_t 
  \dot{\varphi}_t \Omega_t\\
  &  & \\
  & = & \int_X \left( \left\langle \nabla_{g_t}  \dot{\varphi}_t,
  \nabla_{g_t} f_t \right\rangle_{g_t} - 2 f_t  \dot{\varphi}_t 
  \right) \Omega_t\\
  &  & \\
  & = & \int_X \left( \Delta^{\Omega_t}_{g_t} f_t - 2 f_t \right) 
  \dot{\varphi}_t \Omega_t \\
  &  & \\
  & = & \int_X \left( \Delta^{\Omega_t}_{g_t} f_t - 2 f_t +\mathcal{W}_t
  \right)  \dot{\varphi}_t \Omega_t .
\end{eqnarray*}
We set $g \assign g_0$, $\Omega \assign \Omega_0$ and $f \assign f_0$. It is
well known (see for example \cite{Pal2}) that the K\"{a}hler-Ricci soliton condition is
equivalent to the equation 
$$\Delta^{\Omega}_{g_{}} f - 2 f +\mathcal{W}_0 =
0.$$
We deduce
\begin{eqnarray*}
  \ddot{\mathcal{W}}_0 & = & \int_X \frac{d}{d t} \vphantom{dt}_{|_{t = 0}} \left(
  \Delta^{\Omega_t}_{g_t} f_t - 2 f_t +\mathcal{W}_t \right)  \dot{\varphi}_0
  \Omega\\
  &  & \\
  & = & \int_X \frac{d}{d t} \vphantom{dt}_{|_{t = 0}} \left( \Delta^{\Omega_t}_{g_t}
  f_t - 2 f_t \right)  \dot{\varphi}_0 \Omega,
\end{eqnarray*}
thanks to the normalizing condition on $\varphi_t$. We notice now the
following elementary identity obtained in \cite{Pal2} section 3.2
\begin{eqnarray*}
  \frac{d}{d t} \Delta^{\Omega_t}_{g_t} f_t & = & \tmop{div}^{\Omega_t} \left(
  \dot{g}^{\ast}_t \nabla_{g_t} f_t \right) - \left\langle \nabla_{g_t} 
  \dot{\Omega}^{\ast}_t, \nabla_{g_t} f_t \right\rangle_{g_t} +
  \Delta^{\Omega_t}_{g_t}  \dot{f}_t,
\end{eqnarray*}
where $\dot{g}^{\ast}_t : = g^{- 1}_t \dot{g}_t = \omega^{- 1}_t i \partial
\overline{\partial}  \dot{\varphi}_t = \partial^{g_t}_{T_{X, J}} \nabla_{g_t} 
\dot{\varphi}_t$. We infer
\begin{eqnarray*}
  2 \frac{d}{d t} \Delta^{\Omega_t}_{g_t} f_t & = & 2 \tmop{div}^{\Omega_t}
  \left( \partial^{g_t}_{T_{X, J}} \nabla_{g_t}  \dot{\varphi}_t \cdot
  \nabla_{g_t} f_t \right) + 2 \left\langle \nabla_{g_t}  \dot{\varphi}_t,
  \nabla_{g_t} f_t \right\rangle_{g_t} 
  \\
  \\
  &+&\Delta^{\Omega_t}_{g_t}  \left( 2
  \dot{\varphi}_t - \Delta_{g_t}  \dot{\varphi}_t \noplus \right),
\end{eqnarray*}
thanks to (\ref{expVarVol}) and (\ref{evolf}). Applying now the
complex weighted Bochner type identity (13.3) in \cite{Pal2} with $\Omega = \omega^n
/ n!$ we obtain
\begin{eqnarray*}
  2 \partial^{\ast_{g_t}}_{T_{X, J}} \partial^{g_t}_{T_{X, J}} \nabla_{g_t} 
  \dot{\varphi}_t & = & \nabla_{g_t} \Delta_{g_t}  \dot{\varphi}_t,
\end{eqnarray*}
and thus $2 \tmop{div}^{\Omega_t} \partial^{\ast_{g_t}}_{T_{X, J}}
\partial^{g_t}_{T_{X, J}} \nabla_{g_t}  \dot{\varphi}_t = -
\Delta^{\Omega_t}_{g_t} \Delta_{g_t}  \dot{\varphi}_t$. We deduce
\begin{equation}
  \label{varLap} 2 \frac{d}{d t} \Delta^{\Omega_t}_{g_t} f_t = 2
  \tmop{div}^{\Omega_t} \partial^{\ast_{g_t, \Omega_t}}_{T_{X, J}}
  \partial^{g_t}_{T_{X, J}} \nabla_{g_t}  \dot{\varphi}_t + 2 \left\langle
  \nabla_{g_t}  \dot{\varphi}_t, \nabla_{g_t} f_t \right\rangle_{g_t} + 2
  \Delta^{\Omega_t}_{g_t}  \dot{\varphi}_t,
\end{equation}
where $\partial^{\ast_{g_t, \Omega_t}}_{T_{X, J}}$ is the adjoint of
$\partial^{g_t}_{T_{X, J}}$ with respect to the volume form $\Omega_t$. The
proof of the identity (13.3) in \cite{Pal2} shows that at a K\"{a}hler-Ricci soliton
point $\left( J, \omega \right)$ hold the identity
\begin{equation}
  \label{Bochner} 2 \tmop{div}^{\Omega} \partial^{\ast_{g, \Omega}}_{T_{X, J}}
  \partial^g_{T_{X, J}} \nabla_g = - \left( \Delta^{\Omega}_g \right)^2 -
  \left( B^{\Omega}_{g, J} \right)^2,
\end{equation}
where $B^{\Omega}_{g, J} u \assign \tmop{div}^{\Omega} (J \nabla_g u) = g
(\nabla_g u, J \nabla_g f)$. We consider now the complex weighted Laplacian
\begin{eqnarray*}
  \Delta^{\Omega}_{g, J} & = & \Delta^{\Omega}_g - i B^{\Omega}_{g, J} .
\end{eqnarray*}
In the K\"{a}hler-Ricci soliton case hold the identity $[\Delta^{\Omega}_g,
B^{\Omega}_{g, J}] = 0$, (see (15.6) in \cite{Pal2}). This combined with
(\ref{Bochner}) implies
\begin{eqnarray*}
  2 \tmop{div}^{\Omega} \partial^{\ast_{g, \Omega}}_{T_{X, J}}
  \partial^g_{T_{X, J}} \nabla_g & = & - \Delta^{\Omega}_{g, J} 
  \overline{\Delta^{\Omega}_{g, J}} .
\end{eqnarray*}
Using this and the identities (\ref{varLap}) and (\ref{evolf}), we obtain
\begin{eqnarray*}
  \frac{d}{d t} \vphantom{dt}_{|_{t = 0}}  \left( \Delta^{\Omega_t}_{g_t} f_t - 2 f_t
  \right) & = & \left( 2 \Delta_g^{\Omega} + 2 - \frac{1}{2}
  \Delta^{\Omega}_{g, J}  \overline{\Delta^{\Omega}_{g, J}}  \right) 
  \dot{\varphi}_0 \\
  &  & \\
  & = & - \frac{1}{2} P^{\Omega}_{g, J}  \dot{\varphi}_0,
\end{eqnarray*}
where $P^{\Omega}_{g, J} \assign (\Delta^{\Omega}_{g, J} - 2\mathbbm{I})
\overline{ (\Delta^{\Omega}_{g, J} - 2\mathbbm{I})} $ is a non-negative
self-adjoint real elliptic operator with respect to the
$L_{\Omega}^2$-hermitian product (see \cite{Pal2,Pal3}). We deduce
\begin{eqnarray*}
  \ddot{\mathcal{W}}_0 & = & - \frac{1}{2}  \int_X P^{\Omega}_{g, J} 
  \dot{\varphi}_0 \cdot \dot{\varphi}_0 \Omega \leqslant 0 .
\end{eqnarray*}
This is a particular case of Proposition 2 in \cite{Pal2}. The equality hold if and
only if $\dot{\varphi}_0 = \psi + \bar{\psi}$ with $\left( \nabla_{g_0} \psi
\right)_J^{1, 0}$ holomorphic vector field thanks to lemma 22 in \cite{Pal2}.

\begin{remark}
  Notice that the stability result over the space of K\"{a}hler potentials does
  not allow to deduce the general solution in \cite{Pal2,Pal3}. This is because in the non
  K\"{a}hler-Einstein case (it has been pointed out in \cite{Pal2} that in this case the
  solution is trivial), the tangent space of the embedding of the space of
  K\"{a}hler potentials inside the space of complex structures compatible with $\omega$ is not orthogonal and it has positive
  dimensional intersection with the tangent space to the symplectic orbit.
  Furthermore the understanding of the orthogonal behavior of the
  endomorphism Hessian of $\mathcal{W}$ in restriction to such spaces is out
  of reach without using the general solution in \cite{Pal2,Pal3}.
\end{remark}

{\tmstrong{Acknowledgment}}. We invite the readers to compare our note with a
quite recent preprint \cite{Fon}.

\vspace{1cm}
\noindent
Nefton Pali
\\
Universit\'{e} Paris Sud, D\'epartement de Math\'ematiques 
\\
B\^{a}timent 425 F91405 Orsay, France
\\
E-mail: \textit{nefton.pali@math.u-psud.fr}
\end{document}